\documentclass{amsart}

\usepackage{hyperref}

\usepackage[dvips]{graphicx}
\usepackage[latin1]{inputenc}

\input xy
\xyoption{all}

\usepackage{amsmath}
\usepackage{amsfonts}
\usepackage{amssymb}
\usepackage{amsthm}
\usepackage{indentfirst}

\usepackage[T1]{fontenc}

\newcommand{\comment}[1]{}
\DeclareMathOperator\supp{Supp}
\DeclareMathOperator\Tr{Tr}
\DeclareMathOperator\Id{Id}
\DeclareMathOperator\For{For}

\newtheorem{Def}{Definition}[section]
\newtheorem{Defprop}[Def]{Definition-Proposition}
\newtheorem{Th}{Theorem}[section]
\newtheorem{Defth}[Th]{Definition-Theorem}
\newtheorem{Prop}[Th]{Proposition}
\newtheorem{Lem}[Th]{Lemma}
\newtheorem{Corol}[Th]{Corollary}
\newtheorem{Conj}{Conjecture}
\newtheorem{Ques}[Conj]{Question}\theoremstyle{remark}
\newtheorem{Rem}{Remark}

\comment{%version fran?aise
\usepackage[francais]{babel}
\newtheorem{Def}{D?finition}

\newtheorem{Th}{Th?or?me}

\newtheorem{Prop}[Th]{Proposition}
\newtheorem{Lem}[Th]{Lemme}

\theoremstyle{remark}
\newtheorem{Rem}{Remarque}
}

\newcommand{\C}{\mathbb{C}}
\newcommand{\Q}{\mathbb{Q}}

\newcommand{\E}{\mathbb{E}}

\newcommand{\Par}{\mathcal{P}}
\newcommand{\B}{\mathcal{B}}
\newcommand{\A}{\mathcal{A}}

\title[Partial JM elements and star factorizations]{Partial Jucys-Murphy elements and star factorizations}

\author{Valentin F\'eray}
\address{LaBRI, CNRS, Universit\'e Bordeaux 1,
351 cours de la Lib\'eration, 33 400 Talence, France}
\email{feray@labri.fr}


\begin{document}

\maketitle

\begin{abstract}
In this paper, we look at the number of factorizations of a given permutation into star transpositions. In particular, we give a natural explanation of a hidden symmetry, answering a question of I.P. Goulden and D.M. Jackson. We also have a new proof of their explicit formula. Another result is the normalized class expansion of some central elements of the symmetric group algebra introduced by P. Biane.\\

To obtain this results, we use natural analogs of Jucys-Murphy elements in the algebra of partial permutations of V. Ivanov and S. Kerov. We investigate their properties and use a formula of A. Lascoux and J.Y. Thibon to give the expansion of their power sums on the natural basis of the invariant subalgebra.
\end{abstract}

\section{Introduction}

\subsection{Background}

\subsubsection{Counting factorizations}
We denote by $S_n$ the symmetric group over $n$ elements. If $\sigma$ is a permutation and $G$ a set of generators of $S_n$, it is natural to look for the minimal number of factors needed to write $\sigma$ as a product of elements of $G$. This defines a length $\ell_G(\sigma)$ on the symmetric group, which is easy to interpret combinatorially when $G$ is one of the usual sets of generators of $S_n$.\\

A natural question arising is the number of minimal factorizations (decomposition of $\sigma$ as a product of $\ell_G(\sigma)$ generators) : this problem has for instance being studied by J. D\'enes for the set of all transpositions \cite{Denes1959} and by R.P. Stanley with the Coxeter generators \cite{Stanley1984}.\\

One can also look at two interesting variants :
\begin{itemize}
\item counting the number of factorizations of $\sigma$ into $r$ factors in $G$, $r$ being a fixed integer greater than $\ell_G(\sigma)$.
\item counting only transitive factorizations, that is to say factorizations $\sigma = g_1 \ldots g_r$ with the additional condition that the subgroup of $S_n$ spanned by the $g_i$ is the whole symmetric group $S_n$ (this hypothesis might seem strange but it is in fact only the connectedness of a certain ramified covering associated to the factorization).
\end{itemize}

\subsubsection{Transitive star factorization} 
In this article, we look at the set of generators $\{(i\ n),i<n\}$ (called star transpositions) and not limit ourselves to minimal factorizations. Note that, in this context, the transitivity hypothesis consists just in using at least once each generator. This problem has already lead to some interesting papers :
\begin{itemize}
 \item first, I. Pak \cite{Pak1999} has computed the number of minimal factorizations into star transpositions for certain permutations.
 \item second, J. Irving and A. Rattan \cite{IrvingRattan2009} have generalized Pak's result by counting the number of minimal transitive factorizations of any permutation into star transpositions.
 \item third, I. P. Goulden and D. M. Jackson \cite{GouldenJackson2009} have generalized Irving and Rattan's result : they give the number of transitive factorizations of any permutation into a product of any number of star transpositions.
\end{itemize}
As the number of (non necessarily transitive) star factorizations  can be easily recovered from the number of transitive ones, the counting problem is closed for this set of generators.\\

But, in the two papers \cite{IrvingRattan2009,GouldenJackson2009}, the authors notice the following strange fact : although $n$ clearly plays a particular role (as it is in the support of all the generators), the number of transitive star factorizations of a given length depends only on the cycle type of the permutation. In this paper, we provide a natural explanation (without computing the numbers of factorizations) of this phenomenon. This answers a question of these authors.

\subsubsection{Powers of Jucys-Murphy elements}
Star transpositions happen to be very important in representation theory through powers of Jucys-Murphy elements. For instance, A. Okounkov \cite{Okounkov2000} has used the latter to prove the Baik-Deift-Johansson conjecture concerning the fluctuations of the first lines of a random Young diagram under Plancherel's measure. In fact, he only used the combinatorics of decompositions into star factorizations of the identity element, because it is the only one that has a non-zero trace in the regular representation (which corresponds to Plancherel's measure).\\

% P. Biane \cite{Biane1998} have used products of the kind $(\star\ i_1) \ldots (\star\ i_r)$ in $S(\{1,\ldots,n,\star\})$ which fixes the element $\star$ (see paragraph \ref{subsect:moment_transition}) to find the equivalent of the character values in some asymtotic models. 
Instead of powers of Jucys-Murphy elements, one can also look to the following elements introduced by P. Biane \cite{Biane1998} :
$$M_n^r=\E(X_{n+1}^r),$$
where $X_{n+1}$ is the $n+1$-th Jucys-Murphy element in $\Q[S_{n+1}]$ ($X_{n+1}=(1\ n+1) + \ldots +(n\ n+1)$) and $\E$ is the linear map from $\Q[S_{n+1}]$ to $\Q[S_n]$ defined by :
$$\E(\sigma) =
\begin{cases}
 \sigma_{/\{1,\ldots,n\}} & \text{ if } \sigma(n+1)=n+1 ;\\
 0 &\text{ else.}
\end{cases}$$
This element is clearly in the center of the symmetric group algebra $\Q[S_n]$ and it is natural to look at its coefficients on the natural basis $(C_\lambda)_{\lambda \vdash n}$ ($C_\lambda$ is the formal sum of permutations of type $\lambda$). In his paper, P. Biane computes the highest degree part (with respect to a gradation that we do not describe here, but it corresponds to the coefficient of $t^0$ in Theorem \ref{th:class_exp_Mnr}).\\

 Here, we give the complete class expansion of these central elements. As the eigenvalues of these elements on an irreducible representation of $S_n$ can be described combinatorially using the shape of the corresponding diagram (see \cite[Proposition 3.3]{Biane1998}), this formula could be used to generalized Okounkov's result to other measures.

\subsection{Content of this article}

\subsubsection{Results}\label{subsubsect:results}
First, we give a direct proof of the following symmetry property:
\begin{Prop}\label{prop:central}
The number of transitive star factorizations into $r$ factors of a permutation depends only on its conjugacy class in the symmetric group.
\end{Prop}
Then we provide a new proof of Goulden's and Jackson's formula for this number:
\begin{Th}\label{th:GJ_formula}
 The number of transitive star factorization of a permutation of type $\lambda$ into $r$ factors is:
\begin{equation}\frac{r!}{n!} \left( \prod_i \lambda_i \right) \ [t^g] f(t)^{n - 2} \prod_i f(\lambda_i t),\end{equation}
where $g=r-(|\lambda|+\ell(\lambda)-2)$ and $f(t)=2 t^{-1} \sinh(t/2)$.
\end{Th}
Another result is the complete class expansion of the elements $M_n^r$ introduced by P. Biane:
\begin{Th}\label{th:class_exp_Mnr}
If, for $|\lambda| \leq n$, $a_{\lambda;n}$ is the multiple of $C_{\lambda 1^{n-|\lambda|}}$ defined by equation (\ref{eq:a_lambdan}), then:
\begin{equation}M_n^r = \sum_{|\lambda| \leq n} \left[ \frac{r!}{(|\lambda|+1)!} \left( \prod_i \lambda_i \right) \ [t^{r-(|\lambda|+\ell(\lambda))}] f(t)^{|\lambda|} \prod_i f(\lambda_i t) \right] a_{\lambda;n}\end{equation}
\end{Th}
This theorem has also been obtained independently by M. Lassalle \cite{LassalleJM}.

\subsubsection{Tools}
We will use the properties of Jucys-Murphy elements to explain that $n$ does not play a particular role when we count the transitive star factorizations of a permutation. Indeed, a similar phenomenon appears with symmetric functions in Jucys-Murphy elements : although the definition is not symmetric in the integers from $1$ to $n$, they are known to belong to the centre of the symmetric group algebra. In particular, one knows that : $$p_r(X_1,\ldots,X_n) = X_1^r + X_2^r + \cdots + X_n^r \in Z(\Q[S_n]).$$
Roughly speaking, if we take the transitive part, only the last term of the sum has a non-zero contribution (in all other terms, $n$ never appears) so we keep exactly the transitive factorizations of permutations as a product $(i\ n), i<n$. Unfortunately, there is no operator in the symmetric group algebra consisting in taking the transitive part. So we have to work with an upper algebra.\\

The algebra of partial permutations defined by Ivanov and Kerov \cite{IvanovKerov1999} is a good object in this context. So we will define some partial Jucys-Murphy elements in this algebra and investigate their properties.\\

To obtain explicit formulas (Theorems \ref{th:GJ_formula} and \ref{th:class_exp_Mnr}), we need the class expansion of the power sums of Jucys-Murphy elements computed by A. Lascoux and J.Y. Thibon in \cite{LascouxThibon2001}.

\subsection{Outline of the paper}

In section \ref{sect:alg_partiel_perm}, we recall the definition of the algebra of partial permutations by Ivanov and Kerov and study natural analogs of Jucys-Murphy elements.\\

In section \ref{sect:counting_factorizations}, we prove the results claimed in paragraph \ref{subsubsect:results}.

\section{Partial Jucys-Murphy elements}\label{sect:alg_partiel_perm}

\subsection{The algebra of partial permutations}\label{subsect:defs}
In this paragraph, we will describe the algebra of partial permutations, which is a good tool to encode transitivity in the context of star factorizations. This algebra was introduced by S. Kerov and V. Ivanov \cite{IvanovKerov1999}, to understand how the connection coefficients $c_{\lambda,\mu}^\nu$ of symmetric groups behave when one adds parts equal to $1$ to the partitions $\lambda$, $\mu$ and $\nu$. The definitions and properties of paragraphs \ref{subsect:defs} and \ref{subsect:action} (and some others) can be found in their paper.\\

\begin{Def}
 A partial permutation is a pair $(\sigma,d)$, where $d$ is a finite set of integers and $\sigma$ a permutation of $d$.\\

We consider the following multiplication:
$$(\sigma,d) \cdot (\sigma',d') = (\tilde{\sigma} \circ \tilde{\sigma'},d \cup d'),$$
where $\tilde{\sigma}$ (respectively $\tilde{\sigma'}$) is the permutation of $d \cup d'$ which acts like $\sigma$ (respectively $\sigma'$) on $d$ (respectively $d'$) and fixes the other points (in the following, we will often forget the tilde and view a permutation of a set $d_0$ as a permutation of a bigger set $d_1 \supseteq d_0$ without changing the notation).\\

The set of partial permutations endowed with this internal law is a semi-group denoted $\Par_\infty$. If $n$ is a natural number, the set of partial permutations of $\{1,\ldots,n\}$ (which means that $d \subset \{1,\ldots,n\}$) is stable, thus it forms a semi-group $\Par_n$. The corresponding semi-group algebras over the rational field will be denoted by $\B_n:= \Q[\Par_n]$ ($n \leq \infty$).
\end{Def}
% 
% \begin{Rem}
%  We can also see partial permutations as permutations with two kinds of fixed points: those in $d$ and those which are not in $d$.
% \end{Rem}

The main advantage of this tower of algebras, comparing to the symmetric group algebras is that, in addition to the inclusion $B_n \hookrightarrow B_{n+1}$, one has a surjective morphism $\B_{n+1} \rightarrow \B_n$, which consists in keeping only partial permutations for which the set $d$ is included in $\{1,\ldots,n\}$ (the other partial permutations are sent to $0$). In this context, $\B_\infty$ is simply the projective limit of the $\B_n$. The different morphisms are summarized in Figure \ref{fig:morphismesB}: The vertical arrows correspond to morphisms $\For_n$ from $\B_n$ to $\Q[S_n]$ consisting in forgetting $d$.

\begin{figure}[ht]
$$\xymatrix{\ar@<.5ex>@{^{(}->}[r] & \B_{n-1} \ar@{>>}[d] \ar@<.5ex>@{^{(}->}[r]  \ar@<.5ex>@{>>}[l] & \B_n \ar@{>>}[d] \ar@<.5ex>@{^{(}->}[r] \ar@<.5ex>@{>>}[l]& \B_{n+1} \ar@{>>}[d] \ar@<.5ex>@{.}[r] \ar@<.5ex>@{>>}[l] & \B_\infty \ar@{>>}[d] \ar@<.5ex>@{.>>}[l]\\
\ar@{^{(}->}[r] & \Q[S_{n-1}] \ar@{^{(}->}[r] & \Q[S_n] \ar@{^{(}->}[r] & \Q[S_{n+1}] \ar@{^{(}.>}[r] & \Q[S_\infty] } $$
\caption{Algebra morphisms between the $\B_n$'s and the $\Q[S_n]$'s.}
\label{fig:morphismesB}
\end{figure}

\subsection{Action of symmetric groups and invariants}\label{subsect:action}
For $n \leq \infty$, there is a canonical action of the symmetric group $S_n$ on the set of partial permutations $\Par_n$:
$$\tau \cdot (\sigma,d)= (\tau \sigma \tau^{-1}, \tau(d)).$$
The elements of $\B_n$ which are fixed by all permutations form a subalgebra of $\B_n$ which will be denoted by $\mathcal{A}_n$.\\

Note that this subalgebra is smaller than the center of $\B_n$, but its projection on $\Q[S_n]$ is exactly $Z(\Q[S_n])$. The morphisms of Figure \ref{fig:morphismesB} can be restricted to the invariant subalgebras, so we have another family of morphisms represented in Figure \ref{fig:morphismesA}. The algebra $\A_\infty$ is the projective limit of the $\A_n$'s.

\begin{figure}[ht]
$$\xymatrix{\ar@<.5ex>@{^{(}->}[r] & \mathcal{A}_{n-1} \ar@{>>}[d] \ar@<.5ex>@{^{(}->}[r]  \ar@<.5ex>@{>>}[l] & \mathcal{A}_n \ar@{>>}[d] \ar@<.5ex>@{^{(}->}[r] \ar@<.5ex>@{>>}[l]& \mathcal{A}_{n+1} \ar@{>>}[d] \ar@<.5ex>@{.}[r] \ar@<.5ex>@{>>}[l] & \mathcal{A}_\infty \ar@{>>}[d] \ar@<.5ex>@{.>>}[l]\\
\ar@{^{(}->}[r] & Z(\Q[S_{n-1}]) \ar@{^{(}->}[r] & Z(\Q[S_n]) \ar@{^{(}->}[r] & Z(\Q[S_{n+1}]) \ar@{^{(}.>}[r] & Z(\Q[S_\infty]) } $$
\caption{Algebra morphisms between the $\mathcal{A}_n$'s and the $Z(\Q[S_n])$'s.}
\label{fig:morphismesA}
\end{figure}

\paragraph{\textit{Basis of $\mathcal{A}_n$}} If $\lambda$ is a partition, one defines :
\begin{equation}\alpha_{\lambda;n} = \sum_{d \subseteq \{1,\ldots,n\}, |d|=|\lambda| \atop \sigma \in S_d \text{ of type }\lambda} (\sigma,d)\end{equation}
Let us make two little remarks:
\begin{itemize}
\item If $n < |\lambda|$ one has $\alpha_{\lambda,n}=0$.
\item The sequence $(\alpha_{\lambda;n})_{n \geq 1}$ is an element $\alpha_\lambda$ of the projective limit $\mathcal{A}_\infty$ of the $\mathcal{A}_n$'s.
\end{itemize}

\begin{Prop} \mbox{}\\
 For any $n \geq 1$, the family $(\alpha_{\lambda;n})_{|\lambda| \leq n}$ forms a linear basis of $\mathcal{A}_n$.\\
 The family $\alpha_\lambda$ forms a linear basis of $\mathcal{A}_\infty$.\\
\end{Prop}

To formulate our results in the symmetric group algebra and not in the partial permutation algebra, one has to compute the image $a_{\lambda;n}$ of $\alpha_{\lambda;n}$ by the canonical morphism $\For_n$:
\begin{equation}\label{eq:a_lambdan}
a_{\lambda;n}=\binom{n-|\lambda|+m_1(\lambda)}{m_1(\lambda)} C_{\lambda 1^{n-|\lambda|}},
\end{equation}
where $m_1(\lambda)$ denotes the number of parts of $\lambda$ equal to $1$. A simple combinatorial argument shows that if we denote as usual $z_\lambda = \left( \prod_i \lambda_i \right) \cdot \left( \prod_i m_i! \right)$, one has:
$$ a_{\lambda;n}=\frac{n (n-1) \cdots (n- |\lambda|+1)}{z_\lambda} \frac{C_{\lambda 1^{n-|\lambda|}}}{\big|C_{\lambda 1^{n-|\lambda|}}\big|}.$$

\paragraph{\textit{A representation-theorical criterion for an element to belong to $\A_n$}} Kerov and Ivanov have proved that the algebra $\B_n$ is semi-simple. Indeed, one has an isomorphism:
$$
\begin{array}{rcl}
 \B_n & \stackrel{\sim}{\longrightarrow} & \bigoplus\limits_{d \subseteq \{1,\ldots,n\}} \Q[S_d] \\
 x & \mapsto & (\phi_d(x))_{d \subseteq \{1,\ldots,n\}}
\end{array}
\text{ where }
\phi_d(\sigma,d') = \begin{cases}
             \sigma & \text{if } d' \subseteq d ;\\
		0 & \text{else.}
            \end{cases}
$$
Therefore, the irreducible representations of $\B_n$ are indexed by pairs $(d,\lambda)$ where $d$ is a subset of $\{1,\ldots,n\}$ and $\lambda$ is a partition of size $|d|$.\\

One can easily characterize the elements of $\A_n$ by their action on the irreducible representations of $\B_n$.
\begin{Prop}
 An element $x \in \B_n$ is in the invariant algebra $\A_n$ iff
\begin{itemize}
 \item it acts like a homothetic transformation on all irreducible representations of $\B_n$ ;
 \item its eigenvalue on the representation indexed by $(d,\lambda)$ depends only on the size of $d$ and the partition $\lambda$ (and not on the set $d$ itself).
\end{itemize}
 A sequence $(a_n)$ such that $(\forall n, a_n \in \A_n)$ belongs to $A_\infty$ if, for any partition $\lambda$ of size $k$, the eigenvalue of $a_n$ on a representation indexed by a pair $(d,\lambda)$ ($|d|=k$) does not depend on $n \geq k$.
\end{Prop}

\subsection{Partial Jucys-Murphy elements}\label{subsect:part_JM_elemts}

Recall that in the symmetric group algebra $\Q[S_n]$, we call Jucys-Muphy elements the following elements $(X_i)_{1 \leq i \leq n}$:
\begin{equation}X_i = \sum_{j<i} (j\ i).\end{equation}
These elements, introduced independently by A. Jucys \cite{Jucys1974} and G. Murphy \cite{Murphy1981}, have very beautiful properties:
\begin{itemize}
 \item They commute.
 \item Symmetric functions in $X_i$ are exactly the central elements of the symmetric group algebra.
 \item Their action on the Young basis of an irreducible representation of the symmetric group has an easy description.
\end{itemize}

We can very easily define some analog in the algebra of partial permutations $\B_n$:
\begin{equation}\xi_i:= \sum_{j<i} \big((j\ i), \{j,i\}\big).\end{equation}
They have also nice properties:

\begin{Prop}
The partial Jucys-Murphy elements $\xi_1,\ldots,\xi_n$ commute.
\end{Prop}
\begin{proof}
 In fact, we will show a stronger statement: for any $i>1$, the Jucys-Murphy element $\xi_i$ commutes with the image of $\B_{i-1}$.\\

Indeed, let $(\sigma,d)$ be one element of the canonical basis of $\B_{i-1}$ and $j<i$. One has 
\begin{align*}
 (\sigma,d) \cdot ((j \ i),\{j,i\})& = \big(\sigma \cdot (j  \ i),(\{j,i\} \cup d)\big)\\
 ((\sigma(j) \ i),\{\sigma(j),i\})\cdot (\sigma,d) &= \big((\sigma(j) \ i) \cdot \sigma,(\{\sigma(j),i\} \cup d)\big)
\end{align*}
But $\sigma \cdot (j  \ i)=(\sigma(j) \ i) \cdot \sigma$ and one has $\{j,i\} \cup d=\{\sigma(j),i\} \cup d$ because:
\begin{itemize}
 \item either $\sigma(j)=j$.
 \item or both $j$ and $\sigma(j)$ are in $d$.
\end{itemize}
Finally $(\sigma,d) \cdot ((j \ i),\{j,i\})=((\sigma(j) \ i),\{j,i\})\cdot (\sigma,d)$ and since $\sigma$ is a bijection of $\{1,\ldots,i-1\}$, we obtain $(\sigma,d) \cdot \xi_i = \xi_i \cdot (\sigma,d)$ by summing over $j<i$.
\end{proof}

Thanks to this proposition, the evaluation $f(\xi_1,\ldots,\xi_n)$ of a symmetric polynomial $f$ is well-defined.

\begin{Prop}\label{prop:sym_Jucys_in_An}
 If $f$ is a symmetric function, $f(\xi_1,\ldots,\xi_n)$ belongs to $\mathcal{A}_n$. Moreover, the sequence $f_n=f(\xi_1,\ldots,\xi_n)$ is an element $f(\Xi)$ of the projective limit $\mathcal{A}_\infty$.
\end{Prop}

\begin{proof}[First proof.]
It is enough to prove this lemma for any algebraic basis of the ring of symmetric functions. As in the case of classical Jucys-Murphy elements, this property is easy to see if $f$ is an elementary function $e_\ell$ because there exists a closed formula. Indeed, one can prove easily that :
$$e_\ell(\xi_1,\ldots,\xi_n) = \sum_{\sigma \in S_n \atop \kappa(\sigma)=n-\ell} (\sigma, \supp(\sigma)),$$
where $\kappa(\sigma)$ denotes the number of cycles of $\sigma$ (included the fixed points) and $\supp(\sigma)$ the set of non-fixed points. For each $n$, this is an element of $\mathcal{A}_n$ which can be written in the following form :
\begin{equation}e_\ell(\xi_1,\ldots,\xi_n) = \sum_{\lambda \vdash l} \alpha_{\lambda+\mathbf{1};n}\end{equation}
where $\lambda +\mathbf{1}$ is the partition obtained from $\lambda$ by adding $1$ to each part. The second part of the proposition is immediate with this formula.
\end{proof}

\begin{proof}[Second proof.]
 We will use the representation-theorical criterion for an element to belong to $\A_n$. It is easy to see that $\phi_d\big(f(\xi_1,\ldots,\xi_n)\big)$ corresponds to the symmetric function $f$ evaluated on the Jucys-Murphy elements of the symmetric group algebra $\Q[S_d]$. So it acts on the irreducible representation indexed by $(d,\lambda)$ as a homothetic transformation of ratio $f(C)$, where $C$ is the alphabet of the contents of the diagram $\lambda$. This eigenvalue does not depend on $d$, so $f(\xi_1,\ldots,\xi_n)$ belongs to $\A_n$. But it does not depend on $n$ either, so the corresponding sequence is in $\A_\infty$.
\end{proof}

\subsection{A transcendent set of generators for $\A_\infty$}
As mentioned above, when we look at the symmetric group algebras for a fixed $n$, one has the property that the centre is spanned by symmetric functions in Jucys-Murphy elements. One may wonder if there is a highest version of this property in the algebra $\mathcal{A}_\infty$.

\begin{Prop}\label{prop:base_alg_Ainfty}
Let $f_1,f_2,\ldots,f_n,\ldots$ be an algebraic basis of the symmetric function algebra. Then $\alpha_1,f_1(\Xi),\ldots,f_n(\Xi),\ldots$ is an algebraic basis of $\mathcal{A}_\infty$.
\end{Prop}

\begin{proof}
It is enough to prove this proposition for any algebraic basis of the symmetric function algebra. So let us choose the power sums ($\forall i \geq 1, \ f_i=p_i$).\\

If $\lambda=1^{m} 2^{n} 3^{o} \ldots$ is a partition, we denote by $P_\lambda$ the corresponding monomial in $\alpha_1,p_1(\Xi),p_2(\Xi),\ldots$
$$P_\lambda = \alpha_1^{m} p_1(\Xi)^{n} p_2(\Xi)^{o} \ldots$$
One has to prove that the $(P_\lambda)_\lambda$ forms a linear basis of $\mathcal{A}_\infty$. But
\begin{multline*}
 P_\lambda = \sum_{g_1,\ldots,g_{m}} (\Id,\{g_1\}) \cdots (\Id,\{g_{m}\})\\
\cdot \sum_{{h_1 < i_1 \atop \vdots} \atop h_{n} < i_{n}} \big( (h_1\ i_1),\{h_1,i_1\} \big) \cdots \big( (h_{n}\ i_{n}),\{h_{n},i_{n}\} \big)  \\
\cdot \sum_{{j_1, k_1 < \ell_1 \atop \vdots} \atop j_{o}, k_{o} < \ell_{o}} \big( (j_1\ \ell_1)(k_1\ \ell_1),\{j_1,k_1,\ell_1\} \big) \cdots \big( (j_{o}\ \ell_{o})(k_{o}\ \ell_{o}),\{j_{o},k_{o},\ell_{o}\} \big)\\
\ldots 
\end{multline*}
When one expands this expression, every partial permutation $(\sigma,d)$ which appears fulfills $|d| \leq |\lambda|$, with equality if and only if all the corresponding $g,h,i,j,k,l$ are distinct. In this case, the type of the permutation $\sigma$ is necessarily $\lambda$. Moreover each permutation $\sigma$ of type $\lambda$ comes from exactly $\prod_i m_i(\lambda)$ different sequences of $g,h,i,j,k,l$'s. Finally:
$$P_\lambda = \left(\prod_i m_i(\lambda)\right) \cdot \alpha_\lambda + \sum_{|\mu| < |\lambda|} c_\mu \cdot \alpha_\mu,$$
where the $c_\mu$ are non-negative integers. Thus, the $P_\lambda$ are a basis of $\mathcal{A}_\infty$ with a triangular change of basis from the $\alpha_\lambda$.
\end{proof}

\begin{Rem}
An element of $Z(\Q[S_n])$ defines a function on Young diagrams of size $n$. Just take, as image of a partition $\lambda$ the corresponding eigenvalue on the representation indexed by $\lambda$. So if we have a sequence of elements $x_n \in Z(\Q[S_n])$, it defines a function on all the Young diagrams.\\
The elements of $\A_\infty$ are sequences of elements in the $\A_n$, so, thanks to the morphism $\A_n \rightarrow Z(\Q[S_n])$, one can associate in a canonical way to any element $x$ of $\A_\infty$ a function $\varphi_x$ on Young diagrams. Moreover this application is a morphism of algebras. Here are some examples:
\begin{align*}
 \varphi_{\alpha_1} : \lambda & \mapsto |\lambda|;\\
 \varphi_{\alpha_\mu} : \lambda & \mapsto \frac{|\lambda| (|\lambda| - 1 ) \cdots (|\lambda| - |\mu| +1)}{z_\mu} \frac{\chi^\lambda(\mu)}{\chi^\lambda(1)};\\
 \varphi_{f(\Xi)} : \lambda & \mapsto f(C_\lambda),
\end{align*}
where $C_\lambda$ is the multiset of contents of the diagram $\lambda$. In particular, the second example implies that the image of $\A_\infty$ is exactly the algebra $\Lambda^\star$ of polynomial functions on the set of Young diagrams introduced by Kerov and Olshanski \cite{KerovOlshanski1994}. Moreover, our proposition implies that, if $f_1,f_2,\ldots,f_n,\ldots$ is an algebraic basis of the symmetric function algebra, then $$\lambda  \mapsto |\lambda|,\lambda  \mapsto f_1(C_\lambda),\ldots,\lambda  \mapsto f_n(C_\lambda),\ldots$$ is an algebraic basis of $\Lambda^\star$. This description is equivalent to their description as supersymmetric functions in Frobenius coordinates.
\end{Rem}

\begin{Rem}
The algebra $\mathcal{A}_\infty$ is very close to the algebra $\mathcal{K}$ of Farahat and Higman \cite{FaharatHigman1959}: one has an isomorphism of $\C$-algebras
\begin{eqnarray*}
\mathcal{A}_\infty & \stackrel{\sim}{\longrightarrow} & \mathcal{K} ;\\
a_\lambda & \longmapsto & \binom{x-|\lambda|+m_1(\lambda)}{m_1(\lambda)} K(\overline{\lambda}),
\end{eqnarray*}
where $\overline{\lambda}$ is obtained from $\lambda$ by erasing all parts equal to $1$. Moreover, the following diagram is commutative.
$$\xymatrix{ \mathcal{A}_\infty \ar@{<->}^\sim[rr] \ar@{>>}[rd] && \mathcal{K} \ar@{>>}[ld] \\ & \Q[S_n] }$$

So Proposition \ref{prop:base_alg_Ainfty} is equivalent to their main theorem.
\end{Rem}

\section{Star factorizations}\label{sect:counting_factorizations}

\subsection{Transitive powers of Jucys-Murphy elements}

\subsubsection{Central elements}\label{subsect:central}
In this paragraph, we give a natural explanation of the surprising symmetry
property of the number of transitive star factorizations.
First note that, in the algebra of partial permutations,
star factorizations are characterized by the following trivial lemma:
\begin{Lem}
 The factorization $\sigma = (i_1 n) \cdot \ldots \cdot (i_r n)$ is transitive if and only if
$$ \big((i_1 n),\{i_1,n\}\big) \cdot \ldots \cdot \big((i_r n),\{i_r,n\}\big) = (\sigma,\{1,\ldots,n\})$$ in the algebra of partial permutations.
\end{Lem}
Thus one will consider the following operator:
$$\Tr_n : \begin{array}{rcl} \B_n &\longrightarrow& \Q[S_n] ;\\
         		(\sigma,d) & \longmapsto & \begin{cases} \sigma & \text{if }d=\{1,\ldots,n\} ;\\
         		                            0& \text{else.}
         		                           \end{cases}
        \end{array}
$$
The definition of $\mathcal{A}_n$ implies that $\Tr_n(\mathcal{A}_n) \subset Z(\Q[S_n])$.\\

In the context of star factorizations, the operator $\Tr_n$ keeps the transitive part of an expression. So the number of transitive star factorizations into $r$ factors of a permutation $\sigma$ is the coefficient of $\sigma$ in $\Tr_n(\xi_n^r)=\Tr_n(p_r(\xi_1,\ldots,\xi_n))$ (clearly, for $i<n$, $\xi_i^r$ has no transitive part). This element was called transitive power of Jucys-Murphy element by Goulden and Jackson.\\

Proposition \ref{prop:central} is now an immediate corollary of the first part of Proposition \ref{prop:sym_Jucys_in_An}.

\subsubsection{Goulden and Jackson's explicit formula}
In the previous paragraph, we have seen that the expression $\Tr_n(p_r(\xi_1,\ldots,\xi_n))$ for the transitive powers of Jucys-Murphy elements gives a natural explanation of the fact that they belong to the center of the symmetric group algebra. But this formula gives more information: the number of transitive factorizations into $r$ factors of a permutation of a given type is simply one coefficient of $p_r(\Xi)$, written in the basis $(\alpha_\lambda)$. But this expansion can be deduced from a result of A. Lascoux and J.Y. Thibon. We will obtain this way a new proof of Goulden's and Jackson's formula.\\

To do this, we will use the fact that the sequence $\big(p_r(\xi_1,\ldots,\xi_n)\big)_{n \geq 1}$ is in the algebra $\mathcal{A}_\infty$ (second part of Proposition \ref{prop:sym_Jucys_in_An}). In other words, there exists coefficients $g_\lambda$ (not depending on $n$) such that:
\begin{equation}\label{eq:partial_pr_expansion}
 \forall n, \ p_r(\xi_1,\ldots,\xi_n) = \sum_\lambda g_\lambda \alpha_{\lambda;n}
\end{equation}
If one applies the linear operator $\Tr_n$, one obtains :
$$ \Tr_n(p_r(\xi_1,\ldots,\xi_n)) = \sum_{|\lambda|=n} g_\lambda C_\lambda.$$
So, if $|\lambda|=n$, the number of transitive star factorizations into $r$ factors of a given permutation of type $\lambda$ is just $g_\lambda$.\\

But one could also apply to (\ref{eq:partial_pr_expansion}) the morphism $\For_n$ from $\B_n$ to $\Q[S_n]$, one has:
\begin{equation}\label{eq:pr_norm_class_expansion}
\forall n, \ p_r(X_1,\ldots,X_n) = \sum_\lambda g_\lambda a_{\lambda;n}.
\end{equation}
Although different partitions $\lambda$ correspond to proportional elements $a_{\lambda;n}$ (for a given $n$), the coefficients $g_\lambda$ are uniquely determined by (\ref{eq:pr_norm_class_expansion}) (see the proof of \cite[Proposition 7.2]{IvanovKerov1999}). So one can extract the number of transitive star factorizations from the class expansion of the non transitive analog $p_r(X_1,\ldots,X_n)$.\\

But, as mentioned before, A. Lascoux and J.Y. Thibon have computed this expansion. Before stating their result, one has to introduce some notations: for a partition $\lambda=\lambda_1,\ldots,\lambda_\ell=1^{m_1(\lambda)} 2^{m_2(\lambda)} \ldots$, we denote
$$
 \phi_\lambda(t) = \frac{(1-q^{-1})^{|\lambda|-1}}{|\lambda|! z_\lambda} \frac{\prod\limits_i q^{\lambda_i} - 1}{q-1}_{|q=e^t}.
$$
Its coefficients in the exponential series expansion will be denoted $\phi_{\lambda;r}$, \textit{i.e.}
$$\phi_\lambda(t)= \sum_{r \geq 0} \phi_{\lambda;r} \frac{t^r}{r!}.$$

We can now copy equation (37) of \cite{LascouxThibon2001}:
\begin{equation}p_r(X_1,\ldots,X_n) = \sum_{\lambda} \phi_{\lambda;r} z_\lambda a_{\lambda;n}\end{equation}
As (\ref{eq:pr_norm_class_expansion}) entirely determines the $g_\lambda$'s, one has:
\begin{equation}\label{eq:glambda_phi}
g_\lambda = z_\lambda \phi_{\lambda;r} = z_\lambda r!\ [t^r] \phi_\lambda(t)
\end{equation}
But $\phi_\lambda$ can be written as follows:
\begin{align*}
 \phi_\lambda(t) &= (e^{-t/2})^{|\lambda|-1} \cdot \frac{(e^{t/2}-e^{-t/2})^{|\lambda|-1}}{|\lambda|! z_\lambda} \cdot \frac{\prod\limits_i (e^{t/2})^{\lambda_i} (e^{\lambda_i t /2} - e^{-\lambda_i t /2})}{e^{t/2}(e^{t/2} - e^{-t/2})} \\
 &= \frac{1}{|\lambda|! z_\lambda} t^{|\lambda|+\ell(\lambda)-2} f(t)^{|\lambda| - 2} \prod_i \lambda_i f(\lambda_i t), \label{eq:phit_sinh}
\end{align*}
where $f(t)=2 t^{-1} \sinh(t/2)$. So equation (\ref{eq:glambda_phi}) becomes :
\begin{equation}g_\lambda= \frac{r!}{|\lambda|!} \left( \prod_i \lambda_i \right) \ [t^g] f(t)^{|\lambda| - 2} \prod_i f(\lambda_i t),\end{equation}
where $g=r-(|\lambda|+\ell(\lambda)-2)$. If $|\lambda|=n$, the number of transitive star factorizations into $r$ factors of a permutation of type $\lambda$ is $g_\lambda$ and the formula above corresponds to Goulden's and Jackson's formula (\cite[Theorem 1.1]{GouldenJackson2009}).

%\subsubsection{A combinatorial description of eigenvalues}
\subsection{Biane's central elements}\label{subsect:moment_transition}
In the previous paragraph, we have given the coefficients of $p_r(\xi_1,\ldots,\xi_n)$ in the basis $(\alpha_\lambda)_{|\lambda| \leq n}$ and used this expansion to count transitive star factorizations. In this paragraph, we will use the same expansion to give the class expansion of the elements $M_n^r$ proving Theorem \ref{th:class_exp_Mnr}.\\

In fact, the elements $M_n^r$ have a natural analog in the algebra of partial permutations: let us define
\begin{align*}
\tilde{\E}_n & :
\begin{array}{rcl}
 \B_{n+1} & \longrightarrow & \B_n\\
 (\sigma,d) & \longmapsto & \begin{cases}
                             (\sigma_{/d'},d') & \text{if }d = d' \sqcup \{n+1\} \text{ and }\sigma(n+1)=n+1;\\
				0 & \text{else};
                            \end{cases}
\end{array}\\
\mathcal{M}_n^r&=\tilde{\E}_n(\xi_{n+1}^r).
\end{align*}
One can check easily that:
\begin{align*}
\forall i < n+1, \tilde{\E}_n(\xi_i^r) &= 0\\
\tilde{\E}_n(\alpha_{\lambda;n+1}) &= \begin{cases}
                                 \alpha_{\lambda \setminus 1;n} &\text{if $\lambda$ has a part equal to }1 ;\\
				 0 &\text{else.}
                                \end{cases}
\end{align*}
Thus, one has:
\begin{align}
 \mathcal{M}_n^r &= \tilde{\E}_n\big(p_r(\xi_1,\ldots,\xi_{n+1})\big) \nonumber\\
&= \sum_{|\lambda| \leq n} g_{\lambda \cup 1} \alpha_{\lambda;n},
\end{align}
where $g_\lambda$ has been defined and computed in the previous paragraph. An immediate corollary of the last formula is the normalized class expansion of $M_n^r$:
\begin{equation}M_n^r = \sum_{|\lambda| \leq n} g_{\lambda \cup 1} a_{\lambda;n}.\end{equation}
This implies Theorem \ref{th:class_exp_Mnr}. Beware that, for a given partition $\mu \vdash n$, terms corresponding to all partitions $\lambda$ obtained from $\mu$ by erasing parts equal to $1$ contribute to the coefficient of $C_\mu$. But this kind of expression is better for some purposes than the usual class expansion because the coefficients do not depend on $n$.

\bibliographystyle{alpha}
\bibliography{../../courant}

\end{document}